\documentclass[12pt]{article}

\usepackage{amssymb,latexsym,amsmath}

\usepackage{graphicx}

\begin{document}

\newcommand{\bfi}{\bfseries\itshape}

\makeatletter

\@addtoreset{figure}{section}

\def\thefigure{\thesection.\@arabic\c@figure}

\def\fps@figure{h, t}

\@addtoreset{table}{bsection}

\def\thetable{\thesection.\@arabic\c@table}

\def\fps@table{h, t}

\@addtoreset{equation}{section}

\def\theequation{\thesubsection.\arabic{equation}}

\makeatother

\newtheorem{thm}{Theorem}[section]

\newtheorem{prop}[thm]{Proposition}

\newtheorem{lema}[thm]{Lemma}

\newtheorem{cor}[thm]{Corollary}

\newtheorem{defi}[thm]{Definition}

\newtheorem{rk}[thm]{Remark}

\newtheorem{exempl}{Example}[section]

\newenvironment{exemplu}{\begin{exempl}  \em}{\hfill $\square$

\end{exempl}}

\newcommand{\comment}[1]{\par\noindent{\raggedright\texttt{#1}

\par\marginpar{\textsc{Comment}}}}

\newcommand{\todo}[1]{\vspace{5 mm}\par \noindent \marginpar{\textsc{ToDo}}\framebox{\begin{minipage}[c]{0.95 \textwidth}

\tt #1 \end{minipage}}\vspace{5 mm}\par}

\newcommand{\ea}{\mbox{{\bf a}}}

\newcommand{\eu}{\mbox{{\bf u}}}

\newcommand{\ueu}{\underline{\eu}}

\newcommand{\ueo}{\overline{u}}

\newcommand{\oeu}{\overline{\eu}}

\newcommand{\ew}{\mbox{{\bf w}}}

\newcommand{\ef}{\mbox{{\bf f}}}

\newcommand{\eF}{\mbox{{\bf F}}}

\newcommand{\eC}{\mbox{{\bf C}}}

\newcommand{\en}{\mbox{{\bf n}}}

\newcommand{\eT}{\mbox{{\bf T}}}

\newcommand{\eL}{\mbox{{\bf L}}}

\newcommand{\eR}{\mbox{{\bf R}}}

\newcommand{\eV}{\mbox{{\bf V}}}

\newcommand{\eU}{\mbox{{\bf U}}}

\newcommand{\ev}{\mbox{{\bf v}}}

\newcommand{\eve}{\mbox{{\bf e}}}

\newcommand{\uev}{\underline{\ev}}

\newcommand{\eY}{\mbox{{\bf Y}}}

\newcommand{\eK}{\mbox{{\bf K}}}

\newcommand{\eP}{\mbox{{\bf P}}}

\newcommand{\eS}{\mbox{{\bf S}}}

\newcommand{\eJ}{\mbox{{\bf J}}}

\newcommand{\eB}{\mbox{{\bf B}}}

\newcommand{\eH}{\mbox{{\bf H}}}

\newcommand{\leb}{\mathcal{ L}^{n}}

\newcommand{\eI}{\mathcal{ I}}

\newcommand{\eE}{\mathcal{ E}}

\newcommand{\hen}{\mathcal{H}^{n-1}}

\newcommand{\eBV}{\mbox{{\bf BV}}}

\newcommand{\eA}{\mbox{{\bf A}}}

\newcommand{\eSBV}{\mbox{{\bf SBV}}}

\newcommand{\eBD}{\mbox{{\bf BD}}}

\newcommand{\eSBD}{\mbox{{\bf SBD}}}

\newcommand{\ecs}{\mbox{{\bf X}}}

\newcommand{\eg}{\mbox{{\bf g}}}

\newcommand{\paromega}{\partial \Omega}

\newcommand{\gau}{\Gamma_{u}}

\newcommand{\gaf}{\Gamma_{f}}

\newcommand{\sig}{{\bf \sigma}}

\newcommand{\gac}{\Gamma_{\mbox{{\bf c}}}}

\newcommand{\deu}{\dot{\eu}}

\newcommand{\dueu}{\underline{\deu}}

\newcommand{\dev}{\dot{\ev}}

\newcommand{\duev}{\underline{\dev}}

\newcommand{\weak}{\stackrel{w}{\approx}}

\newcommand{\mild}{\stackrel{m}{\approx}}

\newcommand{\strong}{\stackrel{s}{\approx}}

\newcommand{\weakdown}{\rightharpoondown}

\newcommand{\opg}{\stackrel{\mathfrak{g}}{\cdot}}

\newcommand{\opunu}{\stackrel{1}{\cdot}}
\newcommand{\opdoi}{\stackrel{2}{\cdot}}

\newcommand{\opn}{\stackrel{\mathfrak{n}}{\cdot}}
\newcommand{\opx}{\stackrel{x}{\cdot}}

\newcommand{\tr}{\ \mbox{tr}}

\newcommand{\Ad}{\ \mbox{Ad}}

\newcommand{\ad}{\ \mbox{ad}}

\title{Contractible groups and linear dilatation structures}

\author{Marius Buliga \\
\\
Institute of Mathematics, Romanian Academy \\
P.O. BOX 1-764, RO 014700\\
Bucure\c sti, Romania\\
{\footnotesize Marius.Buliga@imar.ro}}

\date{This version:  06.06.2007}

\maketitle

\begin{abstract}
 A dilatation structure on a metric space, is a notion in 
 between a group and a differential structure. 
The basic objects of a dilatation structure are dilatations 
(or contractions). 
The axioms of a dilatation structure set the rules of interaction 
between different dilatations. 

There are  two notions of linearity associated to dilatation structures: 
the linearity of a function between two dilatation 
structures and the linearity of the dilatation structure itself. 

Our main result here is a characterization of contractible groups in terms of 
dilatation structures. To a normed conical group (normed contractible group) 
we can naturally associate a linear dilatation structure. Conversely, any 
linear and strong dilatation structure comes from the dilatation structure 
of a normed contractible group.
\end{abstract}


\paragraph{MSC classes:} 	22E20; 20E36; 20F65; 22A10; 51F99

\newpage

\tableofcontents

\section{Introduction}

Dilatation structures on  metric spaces,  introduced in \cite{buligadil1},  
describe the approximate self-similarity properties of a metric space. A 
dilatation structure is a notion related, but more general, to groups  and  
 differential structures.  

The basic objects of a dilatation structure are dilatations (or contractions). 
The axioms of a dilatation structure set the rules of interaction between
different dilatations. 

A metric space $(X,d)$ which admits a strong dilatation structure (definition 
\ref{defweakstrong}) has a metric tangent space at any point $x \in X$ (theorem 
\ref{thcone}), and any such metric tangent space has an algebraic structure 
of a conical group (theorem \ref{tgene}). Conical groups are particular examples
of contraction groups. The structure of contraction groups is known in some detail, due to Siebert
\cite{siebert}, Wang \cite{wang}, Gl\"{o}ckner and Willis \cite{glockwill}, 
Gl\"{o}ckner \cite{glockner} and references therein.

By a classical result of Siebert \cite{siebert}
proposition 5.4, we can characterize the algebraic structure of the metric
tangent spaces associated to dilatation structures of a certain kind: they are 
Carnot groups, that is simply connected Lie 
groups  whose Lie algebra admits a positive graduation 
(corollary \ref{cortang}). 

Carnot groups appear in many situations, in particular  in relation with
sub-riemannian geometry cf. Bella\"{\i}che   \cite{bell}, groups with 
polynomial growth cf. Gromov \cite{gromovgr}, or Margulis type rigidity  
results cf. Pansu \cite{pansu}. It is part of the author program of research to
show that dilatation structures are natural objects in all these mathematical
subjects. In this respect the corollary \ref{cortang} represents a
generalization  of difficult results in sub-riemannian geometry concerning the structure of the 
metric tangent space at a point of a regular subriemannian manifold. 
 
Linearity is also a property which can be explained with the help of a 
 dilatation structure. In the second section of the paper we explain why 
 linearity can be casted in terms of dilatations. There are in fact  two kinds 
 of linearity: the linearity of a function between 
 two dilatation structures (definition \ref{defgl}) and the linearity of the  
 dilatation structure itself (definition \ref{defilin}). 
 
Our main result here is a characterization of contraction groups in terms of
dilatation structures. To a normed conical group (normed contraction group) we 
can naturally associate a linear dilatation structure (proposition
\ref{pexlin}). Conversely, by theorem  \ref{tdilatlin} any linear and strong 
dilatation structure comes  from the  dilatation structure of a 
normed contraction group.

\section{Linear structure in terms of dilatations}

Linearity is  a basic property related to vector spaces. 
For example, if $\mathbb{V}$ is a real, finite dimensional vector space then a 
transformation $A: \mathbb{V} \rightarrow \mathbb{V}$ is linear if it is a  morphism of 
groups $A: (\mathbb{V}, +)  \rightarrow (\mathbb{V}, +)$  and homogeneous with respect to positive scalars. Furthermore, in a normed vector space we can speak about linear continuous transformations.

A transformation is affine if it is a composition of a translation with a linear transformation.  In this paper 
we shall use the umbrella name "linear" for affine transformations too. 

We try here to explain that linearity property can be entirely phrased in terms of dilatations 
of the vector space $\mathbb{V}$. 

 For the vector space $\mathbb{V}$, the dilatation based 
at $x \in \mathbb{V}$, of coefficient $\varepsilon>0$, is the function 
$$\delta^{x}_{\varepsilon}: \mathbb{V} \rightarrow \mathbb{V} \quad , \quad 
\delta^{x}_{\varepsilon} y = x + \varepsilon (-x+y) \quad . $$
For fixed $x$ the dilatations based at $x$ form a one parameter group which
 contracts any bounded neighbourhood of $x$ to a point, uniformly with respect 
to $x$.

The algebraic structure of $\displaystyle \mathbb{V}$ is encoded in 
dilatations. Indeed,  using dilatations we can recover the operation of addition and multiplication by scalars. 
 
For $\displaystyle x,u,v \in \mathbb{V}$ and $\varepsilon>0$ we define the
following compositions of dilatations: 
\begin{equation}
\Delta_{\varepsilon}^{x}(u,v) = \delta_{\varepsilon^{-1}}^{\delta_{\varepsilon}^{x} u}
 \delta^{x}_{\varepsilon} v \quad , \quad 
\Sigma_{\varepsilon}^{x}(u,v) = \delta_{\varepsilon^{-1}}^{x} \delta_{\varepsilon}^{\delta_{\varepsilon}^{x} u}
 (v) \quad , \quad inv^{x}_{\varepsilon}(u) =  \delta_{\varepsilon^{-1}}^{\delta_{\varepsilon}^{x} u} x \quad . 
\label{operations}
\end{equation}

The meaning of this functions becomes clear if we compute: 
$$\Delta_{\varepsilon}^{x}(u,v) =  x+\varepsilon(-x+u) + (-u+v)  \quad , $$
$$\Sigma_{\varepsilon}^{x}(u,v) =  u+ \varepsilon(-u+x) + (-x+v) \quad ,$$
$$inv^{x}_{\varepsilon}(u) = =x+\varepsilon(-x+u) + (-u+x) \quad .$$
As $\varepsilon \rightarrow 0$ we have the limits: 
$$\lim_{\varepsilon\rightarrow 0} \Delta_{\varepsilon}^{x}(u,v) = \Delta^{x}(u,v) = x+(-u+v) \quad ,$$
$$\lim_{\varepsilon\rightarrow 0} \Sigma_{\varepsilon}^{x}(u,v) = \Sigma^{x}(u,v) = u+(-x+v) \quad ,$$
$$\lim_{\varepsilon\rightarrow 0} inv^{x}_{\varepsilon}(u) = inv^{x}(u) = x-u+x \quad , $$
uniform with respect to $x,u,v$ in bounded sets. 
The function $\displaystyle  \Sigma^{x}(\cdot,\cdot)$ is a group operation, namely the addition operation 
translated such that the neutral element is $x$. Thus, for $x=0$, we recover 
the usual addition   operation. The function 
$\displaystyle inv^{x}(\cdot)$ is the inverse function with respect to addition,
 and $\displaystyle  \Delta^{x}(\cdot,\cdot)$ is the difference function.

Notice that for fixed $x, \varepsilon$ the function $\displaystyle  \Sigma^{x}_{\varepsilon}(\cdot,\cdot)$ is not a 
group operation, first of all because it is not associative. Nevertheless, this function satisfies 
a shifted associativity property, namely (see theorem  \ref{opcollection}) 
$$\Sigma_{\varepsilon}^{x}(\Sigma_{\varepsilon}^{x}(u,v),w) = 
\Sigma_{\varepsilon}^{x}(u, \Sigma_{\varepsilon}^{\delta^{x}_{\varepsilon}u}(v,w)) \quad .$$
Also, the inverse function $\displaystyle inv^{x}_{\varepsilon}$ is not involutive, but shifted involutive (theorem  \ref{opcollection}):  
$$inv_{\varepsilon}^{\delta^{x}_{\varepsilon}u}\left( inv^{x}_{\varepsilon} u\right) = u \quad . $$

Affine continuous  transformations $A:\mathbb{V} \rightarrow \mathbb{V}$ admit the following 
description in terms of dilatations. (We could dispense of continuity hypothesis in this situation, but 
we want to illustrate a general point of view, described further in the paper).  

\begin{prop}
A continuous transformation  $A:\mathbb{V} \rightarrow \mathbb{V}$ is affine if and only if for any 
$\varepsilon \in  (0,1)$, $x,y \in \mathbb{V}$ we have 
\begin{equation}
 A \delta_{\varepsilon}^{x} y \ = \ \delta_{\varepsilon}^{Ax} Ay \quad . 
 \label{eq1proplin} 
 \end{equation}
\label{1proplin}
\end{prop}

The proof is a straightforward consequence of representation formul{\ae} 
(\ref{operations}) for the addition, difference and inverse operations  in terms of dilatations.

\section{Dilatation structures}

We present here a brief introduction into the subject of dilatation 
structures. For more details see Buliga \cite{buligadil1}. The results with 
proofs are new.

\subsection{Notations}

Let $\Gamma$ be  a topological separated commutative group  endowed with a continuous group morphism 
$$\nu : \Gamma \rightarrow (0,+\infty)$$ with $\displaystyle \inf \nu(\Gamma)  =  0$. Here $(0,+\infty)$ is 
taken as a group with multiplication. The neutral element of $\Gamma$ is denoted by $1$. We use the multiplicative notation for the operation in $\Gamma$. 

The morphism $\nu$ defines an invariant topological filter on $\Gamma$ (equivalently, an end). Indeed, 
this is the filter generated by the open sets $\displaystyle \nu^{-1}(0,a)$, $a>0$. From now on 
we shall name this topological filter (end) by "0" and we shall write $\varepsilon \in \Gamma \rightarrow 
0$ for $\nu(\varepsilon)\in (0,+\infty) \rightarrow 0$. 

The set $\displaystyle \Gamma_{1} = \nu^{-1}(0,1] $ is a semigroup. We note $\displaystyle 
\bar{\Gamma}_{1}= \Gamma_{1} \cup \left\{ 0\right\}$
On the set $\displaystyle 
\bar{\Gamma}= \Gamma \cup \left\{ 0\right\}$ we extend the operation on $\Gamma$ by adding the rules  
$00=0$ and $\varepsilon 0 = 0$ for any $\varepsilon \in \Gamma$. This is in agreement with the invariance 
of the end $0$ with respect to translations in $\Gamma$.

The space $(X,d)$ is a complete, locally compact metric space. For any $r>0$  
and any $x \in X$ we denote by $B(x,r)$ the open ball of center $x$ and radius 
$r$ in the metric space $X$.

On the metric space $(X,d)$ we work with the topology (and uniformity) induced 
by the distance. For any $x \in X$ we denote by $\mathcal{V}(x)$ the topological
filter of open neighbourhoods of $x$.

\subsection{Axioms of dilatation structures}
The  axioms of  a dilatation structure $(X,d,\delta)$ are listed further. 
The first axiom is merely a preparation for the next axioms. That is why we 
counted it as axiom 0.

\
\begin{enumerate}
\item[{\bf A0.}] The dilatations $$ \delta_{\varepsilon}^{x}: U(x) 
\rightarrow V_{\varepsilon}(x)$$ are defined for any 
$\displaystyle \varepsilon \in \Gamma, \nu(\varepsilon)\leq 1$. The sets 
$\displaystyle U(x), V_{\varepsilon}(x)$ are open neighbourhoods of $x$.  
All dilatations are homeomorphisms (invertible, continuous, with 
continuous inverse). 

We suppose  that there is a number  $1<A$ such that for any $x \in X$ we have 
$$\bar{B}_{d}(x,A) \subset U(x)  \ .$$
 We suppose that for all $\varepsilon \in \Gamma$, $\nu(\varepsilon) \in 
(0,1)$, we have 
$$ B_{d}(x,\nu(\varepsilon)) \subset \delta_{\varepsilon}^{x} B_{d}(x,A) 
\subset V_{\varepsilon}(x) \subset U(x) \ .$$

 There is a number $B \in (1,A]$ such that  for 
 any $\varepsilon \in \Gamma$ with $\nu(\varepsilon) \in (1,+\infty)$ the 
 associated dilatation  
$$\delta^{x}_{\varepsilon} : W_{\varepsilon}(x) \rightarrow B_{d}(x,B) \ , $$
is injective, invertible on the image. We shall suppose that 
$\displaystyle  W_{\varepsilon}(x) \in \mathcal{V}(x)$, that   
$\displaystyle V_{\varepsilon^{-1}}(x) \subset W_{\varepsilon}(x) $
and that for all $\displaystyle \varepsilon \in \Gamma_{1}$ and 
$\displaystyle u \in U(x)$ we have
$$\delta_{\varepsilon^{-1}}^{x} \ \delta^{x}_{\varepsilon} u \ = \ u \ .$$
\end{enumerate}

We have therefore  the following string of inclusions, for any $\varepsilon \in \Gamma$, $\nu(\varepsilon) \leq 1$, and any $x \in X$:
$$ B_{d}(x,\nu(\varepsilon)) \subset \delta^{x}_{\varepsilon}  B_{d}(x, A) 
\subset V_{\varepsilon}(x) \subset 
W_{\varepsilon^{-1}}(x) \subset \delta_{\varepsilon}^{x}  B_{d}(x, B) \quad . $$

A further technical condition on the sets  $\displaystyle V_{\varepsilon}(x)$ and $\displaystyle W_{\varepsilon}(x)$  will be given just before the axiom A4. (This condition will be counted as part of 
axiom A0.)

\begin{enumerate}
\item[{\bf A1.}]  We  have 
$\displaystyle  \delta^{x}_{\varepsilon} x = x $ for any point $x$. We also have $\displaystyle \delta^{x}_{1} = id$ for any $x \in X$.

Let us define the topological space
$$ dom \, \delta = \left\{ (\varepsilon, x, y) \in \Gamma \times X \times X 
\mbox{ : } \quad \mbox{ if } \nu(\varepsilon) \leq 1 \mbox{ then } y \in U(x) \,
\, , 
\right.$$ 
$$\left. \mbox{  else } y \in W_{\varepsilon}(x) \right\} $$ 
with the topology inherited from the product topology on 
$\Gamma \times X \times X$. Consider also $\displaystyle Cl(dom \, \delta)$, 
the closure of $dom \, \delta$ in $\displaystyle \bar{\Gamma} \times X \times X$ with product topology. 
The function $\displaystyle \delta : dom \, \delta \rightarrow  X$ defined by 
$\displaystyle \delta (\varepsilon,  x, y)  = \delta^{x}_{\varepsilon} y$ is continuous. Moreover, it can be continuously extended to $\displaystyle Cl(dom \, \delta)$ and we have 
$$\lim_{\varepsilon\rightarrow 0} \delta_{\varepsilon}^{x} y \, = \, x \quad . $$

\item[{\bf A2.}] For any  $x, \in K$, $\displaystyle \varepsilon, \mu \in \Gamma_{1}$ and $\displaystyle u \in 
\bar{B}_{d}(x,A)$   we have: 
$$ \delta_{\varepsilon}^{x} \delta_{\mu}^{x} u  = \delta_{\varepsilon \mu}^{x} u  \ .$$

\item[{\bf A3.}]  For any $x$ there is a  function $\displaystyle (u,v) \mapsto d^{x}(u,v)$, defined for any $u,v$ in the closed ball (in distance d) $\displaystyle 
\bar{B}(x,A)$, such that 
$$\lim_{\varepsilon \rightarrow 0} \quad \sup  \left\{  \mid \frac{1}{\varepsilon} d(\delta^{x}_{\varepsilon} u, \delta^{x}_{\varepsilon} v) \ - \ d^{x}(u,v) \mid \mbox{ :  } u,v \in \bar{B}_{d}(x,A)\right\} \ =  \ 0$$
uniformly with respect to $x$ in compact set. 

\end{enumerate}

\begin{rk}
The "distance" $d^{x}$ can be degenerated: there might exist  
$\displaystyle v,w \in U(x)$ such that $\displaystyle d^{x}(v,w) = 0$. 
\label{imprk}
\end{rk}

For  the following axiom to make sense we impose a technical condition on the co-domains $\displaystyle V_{\varepsilon}(x)$: for any compact set $K \subset X$ there are $R=R(K) > 0$ and 
$\displaystyle \varepsilon_{0}= \varepsilon(K) \in (0,1)$  such that  
for all $\displaystyle u,v \in \bar{B}_{d}(x,R)$ and all $\displaystyle \varepsilon \in \Gamma$, $\displaystyle  \nu(\varepsilon) \in (0,\varepsilon_{0})$,  we have 
$$\delta_{\varepsilon}^{x} v \in W_{\varepsilon^{-1}}( \delta^{x}_{\varepsilon}u) \ .$$

With this assumption the following notation makes sense:
$$\Delta^{x}_{\varepsilon}(u,v) = \delta_{\varepsilon^{-1}}^{\delta^{x}_{\varepsilon} u} \delta^{x}_{\varepsilon} v . $$
The next axiom can now be stated: 
\begin{enumerate}
\item[{\bf A4.}] We have the limit 
$$\lim_{\varepsilon \rightarrow 0}  \Delta^{x}_{\varepsilon}(u,v) =  \Delta^{x}(u, v)  $$
uniformly with respect to $x, u, v$ in compact set. 
\end{enumerate}

\begin{defi}
A triple $(X,d,\delta)$ which satisfies A0, A1, A2, A3, but $\displaystyle d^{x}$ is degenerate for some 
$x\in X$, is called degenerate dilatation structure. 

If the triple $(X,d,\delta)$ satisfies A0, A1, A2, A3 and 
 $\displaystyle d^{x}$ is non-degenerate for any $x\in X$, then we call it  a 
 dilatation structure. 

 If a  dilatation structure satisfies A4 then we call it strong dilatation 
 structure. 
 \label{defweakstrong}
\end{defi}

\subsection{Groups with dilatations. Conical groups}

Metric tangent spaces sometimes have a group structure which is compatible 
with dilatations. This structure, of a group with dilatations, is interesting 
by itself. The notion has been introduced in \cite{buliga2}; we describe  it 
further.

  The  following description of local uniform groups 
  is slightly non canonical, but is  motivated by the case of a 
Lie group endowed with a Carnot-Caratheodory  distance induced by a left invariant distribution 
(see for example  \cite{buliga2}). 

We begin with some notations. Let $G$ be a group. 
We introduce first the double of $G$, as the group $G^{(2)} \ = \ G \times G$ 
with operation
$$(x,u) (y,v) \ = \ (xy, y^{-1}uyv) \quad .$$
The operation on the group $G$, seen as the function
$\displaystyle op: G^{(2)} \rightarrow G$, $\displaystyle op(x,y)  =  xy$
is a group morphism. Also the inclusions:
$$i': G \rightarrow G^{(2)} \ , \ \ i'(x) \ = \ (x,e) $$
$$i": G \rightarrow G^{(2)} \ , \ \ i"(x) \ = \ (x,x^{-1}) $$
are group morphisms.

\begin{defi}
\begin{enumerate}
\item[1.]
$G$ is an uniform group if we have two uniformity structures, on $G$ and
$G\times G$,  such that $op$, $i'$, $i"$ are uniformly continuous.

\item[2.] A local action of a uniform group $G$ on a uniform  pointed space $(X, x_{0})$ is a function
$\phi \in W \in \mathcal{V}(e)  \mapsto \hat{\phi}: U_{\phi} \in \mathcal{V}(x_{0}) \rightarrow
V_{\phi}  \in \mathcal{V}(x_{0})$ such that:
\begin{enumerate}
\item[(a)] the map $(\phi, x) \mapsto \hat{\phi}(x)$ is uniformly continuous from $G \times X$ (with product uniformity)
to  $X$,
\item[(b)] for any $\phi, \psi \in G$ there is $D \in \mathcal{V}(x_{0})$
such that for any $x \in D$ $\hat{\phi \psi^{-1}}(x)$ and $\hat{\phi}(\hat{\psi}^{-1}(x))$ make sense and   $\hat{\phi \psi^{-1}}(x) = \hat{\phi}(\hat{\psi}^{-1}(x))$.
\end{enumerate}

\item[3.] Finally, a local group is an uniform space $G$ with an operation defined in a neighbourhood of $(e,e) \subset G \times G$ which satisfies the uniform group axioms locally.
\end{enumerate}
\label{dunifg}
\end{defi}

An uniform group, according to the definition \eqref{dunifg}, is a group $G$ 
such that left translations are uniformly continuous functions and the left 
action of $G$ on itself is uniformly continuous too.

\begin{defi}
A group with dilatations $(G,\delta)$ is a local uniform group $G$  with  a local action of $\Gamma$ (denoted by $\delta$), on $G$ such that
\begin{enumerate}
\item[H0.] the limit  $\displaystyle \lim_{\varepsilon \rightarrow 0} 
\delta_{\varepsilon} x  =  e$ exists and is uniform with respect to $x$ in a compact neighbourhood of the identity $e$.
\item[H1.] the limit
$$\beta(x,y)  =  \lim_{\varepsilon \rightarrow 0} \delta_{\varepsilon}^{-1}
\left((\delta_{\varepsilon}x) (\delta_{\varepsilon}y ) \right)$$
is well defined in a compact neighbourhood of $e$ and the limit is uniform.
\item[H2.] the following relation holds
$$ \lim_{\varepsilon \rightarrow 0} \delta_{\varepsilon}^{-1}
\left( ( \delta_{\varepsilon}x)^{-1}\right)  =  x^{-1}$$
where the limit from the left hand side exists in a neighbourhood of $e$ and is uniform with respect to $x$.
\end{enumerate}
\label{defgwd}
\end{defi}

These axioms are in fact a particular version of the axioms for a dilatation 
structure.

\begin{defi}
A (local) conical group $N$ is a (local) group with a (local) action of
$\Gamma$ by morphisms $\delta_{\varepsilon}$ such that
$\displaystyle \lim_{\varepsilon \rightarrow 0} \delta_{\varepsilon} x \ = \ e$ for any
$x$ in a neighbourhood of the neutral element $e$.
\end{defi}

A conical group is the infinitesimal version of a group with 
dilatations (\cite{buligadil1} proposition 2).

\begin{prop}
Under the hypotheses H0, H1, H2 $(G,\beta, \delta)$ is a conical group, with operation 
$\beta$ and dilatations $\delta$.
\label{here3.4}
\end{prop}

Any group with dilatations has an associated dilatation structure on it.  In a group with dilatations $(G, \delta)$  we define dilatations based in any point $x \in G$ by 
 \begin{equation}
 \delta^{x}_{\varepsilon} u = x \delta_{\varepsilon} ( x^{-1}u)  . 
 \label{dilat}
 \end{equation}
 
\begin{defi} A normed group with dilatations $(G, \delta, \| \cdot \|)$ is a 
group with dilatations  $(G, \delta)$ endowed with a continuous norm  
function $\displaystyle \|\cdot \| : G \rightarrow \mathbb{R}$ which satisfies 
(locally, in a neighbourhood  of the neutral element $e$) the properties: 
 \begin{enumerate}
 \item[(a)] for any $x$ we have $\| x\| \geq 0$; if $\| x\| = 0$ then $x=e$, 
 \item[(b)] for any $x,y$ we have $\|xy\| \leq \|x\| + \|y\|$, 
 \item[(c)] for any $x$ we have $\displaystyle \| x^{-1}\| = \|x\|$, 
 \item[(d)] the limit 
$\displaystyle \lim_{\varepsilon \rightarrow 0} \frac{1}{\nu(\varepsilon)} \| \delta_{\varepsilon} x \| = \| x\|^{N}$ 
 exists, is uniform with respect to $x$ in compact set, 
 \item[(e)] if $\displaystyle \| x\|^{N} = 0$ then $x=e$.
  \end{enumerate}
  \label{dnco}
  \end{defi}
  
  It is easy to see that if $(G, \delta, \| \cdot \|)$ is a normed group with dilatations then $\displaystyle (G, \beta, \delta, \|\cdot\|^{N})$ is a normed conical group. The norm $\displaystyle \|\cdot\|^{N}$ satisfies the 
  stronger form of property (d) definition \ref{dnco}: for any $\varepsilon >0$ 
  $$ \| \delta_{\varepsilon} x \|^{N} = \varepsilon \| x \|^{N} \quad .$$

In a normed group with dilatations we have a natural left invariant distance given by
\begin{equation}
d(x,y) = \| x^{-1}y\| \quad . 
\label{dnormed}
\end{equation}

The following result is theorem 15 \cite{buligadil1}. 

\begin{thm}
Let $(G, \delta, \| \cdot \|)$ be  a locally compact  normed group with dilatations. Then $(G, \delta, d)$ is 
a dilatation structure, where $\delta$ are the dilatations defined by (\ref{dilat}) and the distance $d$ is induced by the norm as in (\ref{dnormed}). 
\label{tgrd}
\end{thm}

\subsection{Carnot groups}

 Normed conical groups generalize the notion of Carnot groups. A simply 
 connected Lie group whose Lie algebra admits a positive graduation is also 
 called a Carnot group. It is in particular
nilpotent. Such objects appear in sub-riemannian
geometry as models of tangent spaces, cf. \cite{bell}, \cite{gromovgr}, \cite{pansu}.

\begin{defi}
A Carnot (or stratified nilpotent) group is a pair $\displaystyle 
(N, V_{1})$ consisting of a real 
connected simply connected group $N$  with  a distinguished subspace  
$V_{1}$ of  the Lie algebra $Lie(N)$, such that  the following   
direct sum decomposition occurs: 
$$n \ = \ \sum_{i=1}^{m} V_{i} \ , \ \ V_{i+1} \ = \ [V_{1},V_{i}]$$
The number $m$ is the step of the group. The number $\displaystyle Q \ = \ \sum_{i=1}^{m} i 
\ dim V_{i}$ is called the homogeneous dimension of the group. 
\label{dccgroup}
\end{defi}

Because the group is nilpotent and simply connected, the
exponential mapping is a diffeomorphism. We shall identify the 
group with the algebra, if is not locally otherwise stated.

The structure that we obtain is a set $N$ endowed with a Lie
bracket and a group multiplication operation, related by the 
Baker-Campbell-Hausdorff formula. Remark that  the group operation
is polynomial.

Any Carnot group admits a one-parameter family of dilatations. For any 
$\varepsilon > 0$, the associated dilatation is: 
$$ x \ = \ \sum_{i=1}^{m} x_{i} \ \mapsto \ \delta_{\varepsilon} x \
= \ \sum_{i=1}^{m} \varepsilon^{i} x_{i}$$
Any such dilatation is a group morphism and a Lie algebra morphism. 

In fact the class of Carnot groups is characterised by the existence of 
dilatations (see Folland-Stein \cite{fostein}, section 1). 

\begin{prop}
Suppose that the Lie algebra $\mathfrak{g}$ admits an one parameter group 
$\varepsilon \in (0,+\infty) \mapsto \delta_{\varepsilon}$ of simultaneously 
diagonalisable Lie algebra isomorphisms. Then $\mathfrak{g}$ is the algebra of 
a Carnot group. 
\end{prop}

We shall construct a norm on a Carnot group $N$. First pick an euclidean norm 
$\| \cdot \|$ on $\displaystyle V_{1}$.  We shall endow the group $N$ with 
a structure of a sub-Riemannian manifold now. For this take the distribution 
obtained from left translates of the space $V_{1}$. The metric on that 
distribution is obtained by left translation of the inner product restricted to
$V_{1}$.

Because $V_{1}$  generates (the algebra) $N$
then any element $x \in N$ can be written as a product of 
elements from $V_{1}$.  An useful lemma is the following (slight reformulation of Lemma 1.40, 
Folland, Stein \cite{fostein}).  

\begin{lema}
Let $N$ be a Carnot group and $X_{1}, ..., X_{p}$ an orthonormal basis 
for $V_{1}$. Then there is a  
 a natural number $M$ and a function $g: \left\{ 1,...,M \right\} 
\rightarrow \left\{ 1,...,p\right\}$ such that any element 
$x \in N$ can be written as: 
\begin{equation}
x \ = \ \prod_{i = 1}^{M} \exp(t_{i}X_{g(i)})
\label{fp2.4}
\end{equation}
Moreover, if $x$ is sufficiently close (in Euclidean norm) to
$0$ then each $t_{i}$ can be chosen such that $\mid t_{i}\mid \leq C 
\| x \|^{1/m}$
\label{p2.4}
\end{lema}

As a consequence we get: 

\begin{cor}
 The Carnot-Carath\'eodory distance 
$$d(x,y) \ = \ \inf \left\{ \int_{0}^{1} \| c^{-1} \dot{c} \| \mbox{ 
d}t \ \mbox{ : } \ c(0) = x , \ c(1) = y , \quad \quad \right.$$ 
$$\left. \quad \quad \quad \quad  c^{-1}(t) \dot{c}(t) \in V_{1} 
\mbox{ for a.e. } t \in [0,1] 
\right\}$$ 
is finite for any two $x,y \in N$.  The distance is obviously left
invariant, thus it induces a norm on $N$. 
\end{cor}

\subsection{Contractible groups}

Conical groups are particular examples of (local) contraction groups. 

\begin{defi}
A contraction group is a pair $(G,\alpha)$, where $G$ is a  
topological group with neutral element denoted by $e$, and $\alpha \in Aut(G)$ 
is an automorphism of $G$ such that: 
\begin{enumerate}
\item[-] $\alpha$ is continuous, with continuous inverse, 
\item[-] for any $x \in G$ we have the limit $\displaystyle 
\lim_{n \rightarrow \infty} \alpha^{n}(x) = e$. 
\end{enumerate}
\end{defi}

We shall be interested in locally compact contraction groups $(G,\alpha)$, such
that $\alpha$ is compactly contractive, that is: for each compact set 
$K \subset G$ and open set $U \subset G$, with $e \in U$, there is 
$\displaystyle N(K,U) \in \mathbb{N}$ such that for any $x \in K$ and $n \in
\mathbb{N}$, $n \geq N(K,U)$, we have $\displaystyle \alpha^{n}(x) \in U$. 
If $G$ is locally compact then a 
necessary and sufficient condition for $(G,\alpha)$ to be compactly contractive
is: $\alpha$ is an uniform contraction, that is each identity neighbourhood 
of $G$ contains an $\alpha$-invariant neighbourhood.

A conical group is an example of a locally compact, compactly contractive,
contraction group. 
Indeed, it suffices to associate to a conical group $(G,\delta)$ the 
contraction group $\displaystyle (G,\delta_{\varepsilon})$, for a fixed 
$\varepsilon \in \Gamma$ with $\nu(\varepsilon) < 1$.

Conversely, to any contraction group $(G,\alpha)$, which is locally compact and 
compactly contractive, associate the conical group $(G,\delta)$, with 
$\displaystyle \Gamma = \left\{ \frac{1}{2^{n}} \mbox{ : } n \in \mathbb{N}
\right\}$ and for any $n \in \mathbb{N}$ and $x \in G$ 
$$\displaystyle \delta_{\frac{1}{2^{n}}} x \ = \ \alpha^{n}(x) \quad . $$

Finally, a local conical
group has only locally the structure of a contraction group. 
The structure of contraction groups is known in some detail, due to Siebert
\cite{siebert}, Wang \cite{wang}, Gl\"{o}ckner and Willis \cite{glockwill}, 
Gl\"{o}ckner \cite{glockner} and references therein. 

For this paper the following results are of interest. We begin with  the
definition of a contracting automorphism group \cite{siebert}, definition 5.1. 

\begin{defi}
Let $G$ be a locally compact group. An automorphism group on $G$ is a family 
$\displaystyle T= \left( \tau_{t}\right)_{t >0}$ in $Aut(G)$, such that 
$\displaystyle \tau_{t} \, \tau_{s} = \tau_{ts}$ for all $t,s > 0$. 

The contraction group of $T$ is defined by 
$$C(T) \ = \ \left\{ x \in G \mbox{ : } \lim_{t \rightarrow 0} \tau_{t}(x) = e
\right\} \quad .$$
The automorphism group $T$ is contractive if $C(T) = G$. 
\end{defi}

It is obvious that a contractive automorphism  group $T$ induces on $G$ a 
structure of conical group. Conversely, any conical group with $\Gamma = 
(0,+\infty)$ has an associated contractive automorphism group (the group of 
dilatations based at the neutral element). 

Further is proposition 5.4 \cite{siebert}. 

\begin{prop}
For a locally compact group $G$ the following assertions are equivalent: 
\begin{enumerate}
\item[(i)] $G$ admits a contractive automorphism group;
\item[(ii)] $G$ is a simply connected Lie group whose Lie algebra admits a 
positive graduation.
\end{enumerate}
\label{psiebert}
\end{prop}

\section{Properties of dilatation structures}

\subsection{First properties}
\label{induced}

The sum, difference, inverse operations induced by a dilatation structure 
give to the space $X$ almost the structure of an affine space. We collect some results from 
\cite{buligadil1} section 4.2 , regarding the properties of these operations.

\begin{thm}
Let $(X,d,\delta)$ be a  dilatation structure. Then, for any $x \in X$, $\varepsilon \in \Gamma$, 
$\nu(\varepsilon) < 1$,     we have: 
\begin{enumerate}
\item[(a)] for any $u\in U(x)$,  $\displaystyle \Sigma_{\varepsilon}^{x}(x,u) = u $  .
\item[(b)] for any $u\in U(x)$ the functions $\displaystyle  \Sigma_{\varepsilon}^{x}(u,\cdot)$ and 
$\displaystyle  \Delta_{\varepsilon}^{x}(u, \cdot)$ are inverse one to another. 
\item[(c)] the inverse function is shifted involutive:  for any $u\in U(x)$, 
$$inv^{\delta_{\varepsilon}^{x} u}_{\varepsilon} \, inv^{x}_{\varepsilon} (u) = u \quad . $$
\item[(d)] the sum operation is shifted associative: for  any $u, v, w$ sufficiently close to $x$ we have 
$$\Sigma_{\varepsilon}^{x} \left( u, \Sigma_{\varepsilon}^{\delta^{x}_{\varepsilon} u} (v, w)\right) = 
\Sigma_{\varepsilon}^{x} ( \Sigma^{x}(u,v), w) \quad . $$
\item[(e)] the difference, inverse and sum operations are related by 
$$ \Delta_{\varepsilon}^{x}(u,v) = \Sigma_{\varepsilon}^{\delta_{\varepsilon}^{x} u}
 \left( inv_{\varepsilon}^{x}(u), v\right) \quad , $$
 for any $u,v$ sufficiently close to $x$. 
 \item[(f)] for any $u,v$ sufficiently close to $x$ and $\mu \in \Gamma$, 
$\nu(\mu) < 1$,     we have: 
$$\Delta^{x}_{\varepsilon} \left( \delta^{x}_{\mu} u,  \delta^{x}_{\mu} v \right) =  \delta^{ \delta^{x}_{\epsilon \mu} u}_{\mu} \Delta^{x}_{\varepsilon \mu} (u,v) \quad . $$
\end{enumerate}
\label{opcollection}
\end{thm}

\subsection{Tangent bundle}

A reformulation of parts of theorems 6,7 \cite{buligadil1} is the following. 

\begin{thm}
A   dilatation structure $(X,d,\delta)$ has the following properties.  
\begin{enumerate}
\item[(a)] For all $x\in X$, $u,v \in X$ such that $\displaystyle d(x,u)\leq 1$ and $\displaystyle d(x,v) \leq 1$  and all $\mu \in (0,A)$ we have: 
$$d^{x}(u,v) \ = \ \frac{1}{\mu} d^{x}(\delta_{\mu}^{x} u , \delta^{x}_{\mu} v) \ .$$
We shall say that $d^{x}$ has the cone property with respect to dilatations. 
\item[(b)] The metric space $(X,d)$ admits a metric tangent space 
at $x$, for any point $x\in X$. More precisely we have  the following limit: 
$$\lim_{\varepsilon \rightarrow 0} \ \frac{1}{\varepsilon} \sup \left\{  \mid d(u,v) - d^{x}(u,v) \mid \mbox{ : } d(x,u) \leq \varepsilon \ , \ d(x,v) \leq \varepsilon \right\} \ = \ 0 \ .$$
\end{enumerate}
\label{thcone}
\end{thm}

For the next theorem (composite of results in theorems 8, 10 \cite{buligadil1}) 
 we need the previously introduced notion  of a normed conical (local) group.

\begin{thm}
Let $(X,d,\delta)$ be a strong dilatation structure. Then for any $x \in X$ the triple 
 $\displaystyle (U(x), \Sigma^{x}, \delta^{x})$ is a normed local conical group,
 with the norm induced by the distance $\displaystyle d^{x}$.
\label{tgene}
\end{thm}

The conical group $\displaystyle (U(x), \Sigma^{x}, \delta^{x})$ can be regarded as the tangent space 
of $(X,d, \delta)$ at $x$. Further will be denoted by: 
$\displaystyle T_{x} X =  (U(x), \Sigma^{x}, \delta^{x})$.

\begin{defi}
Let $(X,\delta,d)$ be a   dilatation structure and $x\in X$ a point. 
In a neighbourhood $U(x)$ of $x$, for  any $\mu\in (0,1)$ 
we defined the distances:
$$(\delta^{x},\mu)(u,v) = \frac{1}{\mu} d(\delta^{x}_{\mu} u , \delta^{x}_{\mu} v) . $$
\label{drelative}
\end{defi}

\begin{prop}
Let $(X,\delta,d)$ be a  (strong) dilatation structure. 
For any $u,v\in U(x)$ let us  define 
$$\hat{\delta}_{\varepsilon}^{u} v = \Sigma^{x}_{\mu}(u, \delta^{\delta^{x}_{\mu}u}_{\varepsilon} 
\Delta^{x}_{\mu}(u,v)) = \delta^{x}_{\mu^{-1}} \delta_{\varepsilon}^{\delta_{\mu}^{x} u} \delta_{\mu}^{x} v  . $$
Then $\displaystyle (U(x),\hat{\delta}, (\delta^{x},\mu))$ is a  (strong) dilatation structure. 
\label{pshift}
\end{prop}

\paragraph{Proof.}
We have to check the axioms.  The first part of axiom A0 is an easy consequence of theorem 
\ref{thcone} for $(X,\delta,d)$. The second part of A0, A1 and A2 are true based on simple computations. 

The first interesting fact is related to axiom A3. Let us compute, for $v,w \in U(x)$, 
$$\frac{1}{\varepsilon} (\delta^{x},\mu)(\hat{\delta}^{u}_{\varepsilon} v, \hat{\delta}^{u}_{\varepsilon} w) = 
\frac{1}{\varepsilon \mu} d( \delta^{x}_{\mu} \hat{\delta}^{u}_{\varepsilon} v, \delta^{x}_{\mu}
\hat{\delta}^{u}_{\varepsilon} w) = $$
$$ = \frac{1}{\varepsilon \mu} d( \delta_{\varepsilon}^{\delta^{x}_{\mu} u} \delta_{\mu}^{x} v ,  
\delta_{\varepsilon}^{\delta^{x}_{\mu} u} \delta_{\mu}^{x} w) = \frac{1}{\varepsilon\mu} d(  
\delta_{\varepsilon \mu}^{\delta^{x}_{\mu} u} \Delta_{\mu}^{x}(u,v),  \delta_{\varepsilon \mu}^{\delta^{x}_{\mu} u} \Delta_{\mu}^{x}(u,w)) =  $$
$$= (\delta^{\delta_{\mu}^{x} u}, \varepsilon \mu) (  \Delta_{\mu}^{x}(u,v) ,  \Delta_{\mu}^{x}(u,w)) . $$
The axiom A3 is then a consequence of axiom A3 for $(X,\delta,d)$ and we have 
$$\lim_{\varepsilon\rightarrow 0} \frac{1}{\varepsilon} (\delta^{x},\mu)(\hat{\delta}^{u}_{\varepsilon} v, \hat{\delta}^{u}_{\varepsilon} w) = d^{\delta_{\mu}^{x} u} (  \Delta_{\mu}^{x}(u,v) ,  \Delta_{\mu}^{x}(u,w)) . $$
The axiom A4 is also a straightforward consequence of A4 
for $(X,\delta,d)$ and is left to the reader. \quad $\square$

 The proof of the following proposition is an easy computation, of the same 
type as in the lines above, therefore we shall not write it here. 

\begin{prop}
With the same notations as in proposition \ref{pshift}, the transformation 
$\displaystyle \Sigma^{x}_{\mu}(u, \cdot)$ is an isometry from 
$\displaystyle (\delta^{\delta^{x}_{\mu} u}, \mu)$ to 
$\displaystyle (\delta^{x}, \mu)$. Moreover, we have 
$\displaystyle  \Sigma^{x}_{\mu} (u, \delta_{\mu}^{x} u) = u$. 
\label{isoshift}
\end{prop}

These two propositions show  that on a dilatation structure we almost have 
translations 
(the operators $\displaystyle \Sigma^{x}_{\varepsilon}(u, \cdot)$), which are 
almost isometries (that is, not with respect to the distance $d$, but with respect to distances of type $\displaystyle (\delta^{x},\mu)$). 
It is almost as if we were working with  a normed 
conical group, only that we have to use families of distances and to make 
small shifts in the tangent space (as in the last formula in the proof 
of proposition \ref{pshift}).

\subsection{Topological considerations}
 
 In this subsection we compare various topologies and uniformities related to a dilatation structure.  
  
 The axiom A3 implies that for any $x \in X$ the function $\displaystyle d^{x}$ is continuous, therefore 
 open sets with respect to $\displaystyle d^{x}$ are open with respect to $d$. 
 
 If $(X,d)$ is separable and $\displaystyle d^{x}$ is non degenerate then $\displaystyle (U(x), d^{x})$ is also separable and the topologies of $d$ and $\displaystyle d^{x}$ are the same. Therefore $\displaystyle (U(x), d^{x})$ is also locally compact (and a set is compact  with respect to $d^{x}$ if and only if it is compact with respect to $d$). 
 
 If  $(X,d)$ is separable and $\displaystyle d^{x}$ is non degenerate then the uniformities induced by 
 $d$ and  $\displaystyle d^{x}$ are the same. Indeed, let 
 $\displaystyle \left\{u_{n} \mbox{ : } n \in \mathbb{N}\right\}$ 
 be a dense set in $U(x)$, with $\displaystyle x_{0}=x$. 
 We can embed $\displaystyle (U(x), 
 (\delta^{x}, \varepsilon))$ isometrically in the separable Banach space 
 $\displaystyle l^{\infty}$, for any $\varepsilon \in (0,1)$, by the function 
 $$\phi_{\varepsilon}(u) = \left( \frac{1}{\varepsilon} d(\delta^{x}_{\varepsilon}u,  \delta^{x}_{\varepsilon}x_{n}) - \frac{1}{\varepsilon} d(\delta^{x}_{\varepsilon}x,  \delta^{x}_{\varepsilon}x_{n})\right)_{n}  . $$
 A reformulation of point (a) in theorem \ref{thcone} is that on compact sets $\displaystyle \phi_{\varepsilon}$ uniformly converges to the isometric embedding of $\displaystyle (U(x), d^{x})$ 
 $$\phi(u) = \left(  d^{x}(u,  x_{n}) - d^{x}(x, x_{n})\right)_{n}  . $$
Remark that the uniformity induced by $(\delta,\varepsilon)$ is the same as the uniformity 
induced by $d$, and that it is the same induced from the uniformity on 
$\displaystyle l^{\infty}$ by 
the embedding $\displaystyle \phi_{\varepsilon}$. We proved that the uniformities induced by 
 $d$ and  $\displaystyle d^{x}$ are the same.

From previous considerations we deduce the following characterisation of  tangent
spaces asociated to a dilatation structure. 

\begin{cor}
Let $(X,d,\delta)$ be a strong dilatation structure with group $\Gamma = (0,+\infty)$. 
Then for any $x \in X$ the local group  
 $\displaystyle (U(x), \Sigma^{x})$ is locally a simply connected Lie group 
 whose Lie algebra admits a positive graduation (a Carnot group).
 \label{cortang}
\end{cor}

\paragraph{Proof.}
Use the facts:  $\displaystyle (U(x), \Sigma^{x})$ is a locally compact
group (from previous topological considerations)  which admits 
$\displaystyle \delta^{x}$ as a contractive automorphism 
group (from theorem \ref{tgene}). Apply then  Siebert proposition \ref{psiebert}.
 \quad $\square$

\section{Linearity and dilatation structures}

\begin{defi}
Let $(X,d,\delta)$ be a    dilatation structure. A transformation $A:X\rightarrow X$ is linear if it 
is Lipschitz and it commutes with dilatations in the following sense: for any $x\in X$, $u \in U(x)$ and 
$\varepsilon \in \Gamma$, $\nu(\varepsilon) < 1$, if  $A(u) \in U(A(x))$ then  
$$ A \delta^{x}_{\varepsilon} = \delta^{A(x)} A(u) \quad  .$$
The  group of linear transformations, denoted by $GL(X,d,\delta)$ is formed by all invertible and 
bi-lipschitz linear transformations of $X$. 
\label{defgl}
\end{defi}

$GL(X,d,\delta)$ is  indeed a (local) group. In order to see this 
 we start from the  
remark that if $A$ is Lipschitz then there exists $C>0$ such that for all $x\in X$ and $u \in B(x,C)$ we have $A(u)\in U(A(x))$.  The inverse of $A \in GL(X,d,\delta)$ is then linear. Same considerations apply for the composition of two linear, bi-lipschitz and invertible transformations. 

In the particular case of  the first subsection of this paper, namely  
$X$ finite dimensional real, normed vector space, 
$d$ the distance given by the norm, $\Gamma = (0,+\infty)$ and dilatations 
$\displaystyle \delta_{\varepsilon}^{x} u = x + \varepsilon(u-x)$, 
a linear transformations in the sense of definition \ref{defgl} is an affine transformation of the vector 
space $X$.

\begin{prop}
Let $(X,d,\delta)$ be a   dilatation structure and $A:X\rightarrow X$ a linear transformation. Then: 
\begin{enumerate}
\item[(a)] for all $x\in X$, $u,v \in U(x)$ sufficiently close to $x$, we have: 
$$A \, \Sigma_{\varepsilon}^{x}(u,v)  = \Sigma_{\varepsilon}^{A(x)}(A(u),A(v)) \quad . $$
\item[(b)] for all $x\in X$, $u \in U(x)$ sufficiently close to $x$, we have: 
$$ A \, inv^{x}(u) = \, inv^{A(x)} A(u) \quad . $$
\end{enumerate}
\label{plinear}
\end{prop}

\paragraph{Proof.}
Straightforward, just use the commutation with dilatations. $\quad \square$

\subsection{Differentiability of linear transformations}

In this subsection we briefly recall the notion of differentiability associated 
to dilatation structures (section 7.2 \cite{buligadil1}). 
Then we apply it for linear transformations. 

First we need the natural definition below. 

\begin{defi}
 Let $(N,\delta)$ and $(M,\bar{\delta})$ be two conical groups. A function $f:N\rightarrow M$ is a conical group morphism if $f$ is a group morphism and for any $\varepsilon>0$ and $u\in N$ we have 
 $\displaystyle f(\delta_{\varepsilon} u) = \bar{\delta}_{\varepsilon} f(u)$. 
\label{defmorph}
\end{defi}

The definition of derivative with respect to dilatations structures follows. 

 \begin{defi}
 Let $(X, \delta , d)$ and $(Y, \overline{\delta} , \overline{d})$ be two 
 strong dilatation structures and $f:X \rightarrow Y$ be a continuous function. The function $f$ is differentiable in $x$ if there exists a 
 conical group morphism  $\displaystyle Q^{x}:T_{x}X\rightarrow T_{f(x)}Y$, defined on a neighbourhood of $x$ with values in  a neighbourhood  of $f(x)$ such that 
\begin{equation}
\lim_{\varepsilon \rightarrow 0} \sup \left\{  \frac{1}{\varepsilon} \overline{d} \left( f\left( \delta^{x}_{\varepsilon} u\right) ,  \overline{\delta}^{f(x)}_{\varepsilon} Q^{x} (u) \right) \mbox{ : } d(x,u) \leq \varepsilon \right\}Ê  = 0 , 
\label{edefdif}
\end{equation}
The morphism $\displaystyle Q^{x}$ is called the derivative of $f$ at $x$ and will be sometimes denoted by $Df(x)$.

The function $f$ is uniformly differentiable if it is differentiable everywhere and the limit in (\ref{edefdif}) 
is uniform in $x$ in compact sets. 
\label{defdif}
\end{defi}

The following proposition has then a straightforward proof. 

\begin{prop}
Let $(X,d,\delta)$ be a strong  dilatation structure and $A:X\rightarrow X$ a linear
transformation. Then  $A$ is uniformly differentiable and the derivative equals $A$.
\end{prop}

\subsection{Linearity of strong dilatation structures}

Remark that for general dilatation structures the "translations" 
$\displaystyle \Sigma^{x}_{\varepsilon}(u, \cdot)$  are not linear. 
Nevertheless, they commute with dilatation in a known way, according 
to point (f) theorem  \ref{opcollection}. This is important, because the 
transformations $\displaystyle \Sigma^{x}_{\varepsilon}(u, \cdot)$ really behave as 
translations, as explained in subsection \ref{induced}.

The reason for which translations are not linear is that dilatations are
generally not linear. Before giving the next definition we need to establish
a simple estimate. Let $K \subset X$ be compact, non empty set. Then there
is a constant $C(K) > 0$, depending on the set $K$ such that for any $\varepsilon,\mu \in 
\Gamma$ with $\nu(\varepsilon),\nu(\mu) \in (0,1]$ and any $x,y,z \in K$ with 
$d(x,y), d(x,z), d(y,z) \leq C(K)$  we have 
$$\delta^{y}_{\mu}z \in V_{\varepsilon}(x) \quad , \quad
\delta_{\varepsilon}^{x} z \in V_{\mu}(\delta^{x}_{\varepsilon} y) \quad .$$
Indeed, this is coming from the uniform (with respect to K) estimates:  
$$d(\delta^{x}_{\varepsilon} y, \delta^{x}_{\varepsilon} z) \leq \varepsilon 
d^{x}(y,z) + \varepsilon \mathcal{O}(\varepsilon) \quad , $$
$$d(x, \delta^{y}_{\mu} z) \leq d(x,y) + d(y, \delta^{y}_{\mu}z) 
\leq d(x,y) + \mu d^{y}(y,z) + \mu \mathcal{O}(\mu) \quad . $$

\begin{defi} 
A property $\displaystyle \mathcal{P}(x_{1},x_{2},
x_{3}, ...)$ holds for $\displaystyle x_{1}, x_{2}, x_{3},
...$ sufficiently closed if for any compact, non empty set $K \subset X$, there
is a positive constant $C(K)> 0$ such that $\displaystyle \mathcal{P}(x_{1},x_{2},
x_{3}, ...)$ is true for any $\displaystyle x_{1},x_{2},
x_{3}, ... \in K$ with $\displaystyle d(x_{i}, x_{j}) \leq C(K)$.  
\end{defi}

For example, we may say that the expressions 
$$\delta_{\varepsilon}^{x} \delta^{y}_{\mu} z \quad , \quad
\delta^{\delta^{x}_{\varepsilon} y}_{\mu} \delta^{x}_{\varepsilon} z$$ 
are well defined for any $x,y,z \in X$ sufficiently closed and for any 
$\varepsilon,\mu \in 
\Gamma$ with $\nu(\varepsilon),\nu(\mu) \in (0,1]$. 

\begin{defi}
A   dilatation structure $(X,d,\delta)$ is linear if for any $\varepsilon,\mu \in 
\Gamma$ such that $\nu(\varepsilon),\nu(\mu) \in (0,1]$, and for any 
$x,y,z \in X$ sufficiently closed we have 
$$\delta_{\varepsilon}^{x} \,  \delta^{y}_{\mu} z \ = \ 
\delta^{\delta^{x}_{\varepsilon} y}_{\mu} \delta^{x}_{\varepsilon} z \quad .$$
\label{defilin}
\end{defi}

Linear dilatation structures are very particular dilatation structures. The next
proposition gives a family of examples of linear dilatation structures. 

\begin{prop}
The dilatation structure associated to a normed conical group is linear.
\label{pexlin}
\end{prop}

\paragraph{Proof.}
Indeed, for the dilatation structure associated to a normed conical group we
have, with the notations from definition \ref{defilin}: 
$$\delta^{\delta^{x}_{\varepsilon} y}_{\mu}
\delta^{x}_{\varepsilon} z \ = \ \left(x \delta_{\varepsilon} (x^{-1}y) \right) 
\, \delta_{\mu} \left( \delta_{\varepsilon} (y^{-1}x) \, x^{-1} \, x \, 
\delta_{\varepsilon} (x^{-1} z) \right) \ =  $$
$$= \  \left(x \delta_{\varepsilon} (x^{-1}y) \right) 
\, \delta_{\mu} \left( \delta_{\varepsilon} (y^{-1}x)  \, 
\delta_{\varepsilon} (x^{-1} z) \right) \ = \ 
\left(x \delta_{\varepsilon} (x^{-1}y) \right) 
\, \delta_{\mu} \left( \delta_{\varepsilon} (y^{-1}z) \right) \ =  $$
$$= \ x \left(\delta_{\varepsilon} (x^{-1}y) 
\, \delta_{\varepsilon} \, \delta_{\mu} (y^{-1}z) \right) \ = \ 
x \, \delta_{\varepsilon} \left(x^{-1}y 
\,  \delta_{\mu} (y^{-1}z)\right)  \ = \ \delta_{\varepsilon}^{x} \, 
\delta^{y}_{\mu}z \quad .$$
Therefore the dilatation structure is linear. \quad $\square$

In the proposition below we give a relation, true for linear dilatation
structures, with an interesting  interpretation. Let us think in affine 
terms: for closed points $x,u,v$, we think about let us denote 
$\displaystyle w = \Sigma^{x}_{\varepsilon}(u,v)$. We may think that 
the "vector" $(x,w)$ is (approximatively, due to the parameter 
$\varepsilon$)  the sum of the vectors $(x,u)$ and $(x,v)$, based at $x$.  
Denote also  $\displaystyle w' = 
\Delta^{u}_{\varepsilon}(x,v)$; then the "vector" $\displaystyle (u,w')$ is 
 (approximatively) equal to  the differrence 
between the vectors $(u,v)$ and $(u,x)$, based at $u$. In a classical  affine 
space we would have $\displaystyle w = w'$. The same is true for a linear dilatation
structure. 

\begin{prop}
For a linear    dilatation structure $(X,\delta,d)$, for any $x,u,v \in X$ 
sufficiently closed and for any $\varepsilon \in \Gamma$, $\nu(\varepsilon) \leq
1$, we have: 
$$\Sigma^{x}_{\varepsilon}(u,v) \ = \ \Delta^{u}_{\varepsilon}(x,v) \quad . $$
\end{prop}

\paragraph{Proof.}
We have the following string of equalities, by using twice the linearity of the dilatation
structure: 
$$\Sigma_{\varepsilon}^{x}(u,v) \ = \ \delta_{\varepsilon^{-1}}^{x}
\delta_{\varepsilon}^{\delta_{\varepsilon}^{x} u} v \ = \
\delta_{\varepsilon}^{u} \, \delta_{\varepsilon^{-1}}^{x} v \ = \ $$
$$= \ \delta_{\varepsilon^{-1}}^{\delta_{\varepsilon}^{u} x}
\delta_{\varepsilon}^{u} v \ = \ \Delta_{\varepsilon}^{u}(x,v) \quad . $$
The proof is done. \quad $\square$

The following expression: 
\begin{equation}
Lin(x,y,z; \varepsilon, \mu) \ = \ d\left( \delta_{\varepsilon}^{x} \, 
\delta^{y}_{\mu}z \, , \,  \delta^{\delta^{x}_{\varepsilon} y}_{\mu}
\delta^{x}_{\varepsilon} z\right) 
\label{deflinfunct}
\end{equation}
is a measure of lack of linearity, for a general    dilatation structure. 
The next theorem shows that, infinitesimally, any  dilatation
structure is linear. 

\begin{thm}
Let $(X,d,\delta)$ be a  strong dilatation structure. Then for any $x, y, z \in X$ 
sufficiently close we have 
\begin{equation}
\lim_{\varepsilon \rightarrow 0} \frac{1}{\varepsilon^{2}} \, 
Lin(x,\delta_{\varepsilon}^{x} y, \delta_{\varepsilon}^{x} z ; \varepsilon,
\varepsilon) \ = \ 0 \quad . 
\label{inflin}
\end{equation}
\label{thinflin}
\end{thm}

\paragraph{Proof.}
From the hypothesis of the theorem we have:
$$\frac{1}{\varepsilon^{2}} \, 
Lin(x,\delta_{\varepsilon}^{x} y, \delta_{\varepsilon}^{x} z ; \varepsilon,
\varepsilon) \ = \ \frac{1}{\varepsilon^{2}} \, d\left( 
\delta^{x}_{\varepsilon} \, \delta^{\delta^{x}_{\varepsilon} y}_{\varepsilon} z 
\, , \, 
\delta_{\varepsilon}^{\delta_{\varepsilon^{2}}^{x} y} \delta_{\varepsilon}^{x} z
\right) \ = \ $$
$$= \ \frac{1}{\varepsilon^{2}} \, d\left( \delta^{x}_{\varepsilon^{2}} \, 
\Sigma^{x}_{\varepsilon}(y,z) \, , \, \delta^{x}_{\varepsilon^{2}} \, 
\delta^{x}_{\varepsilon^{-2}} \, 
\delta_{\varepsilon}^{\delta_{\varepsilon^{2}}^{x} y} \delta_{\varepsilon}^{x} z
\right) \ = \ $$ 
$$= \ \frac{1}{\varepsilon^{2}} \, d\left( \delta^{x}_{\varepsilon^{2}}
\, \Sigma^{x}_{\varepsilon}(y,z) \, , \, \delta^{x}_{\varepsilon^{2}} \, 
\Sigma^{x}_{\varepsilon^{2}} (y \, , \, 
\Delta^{x}_{\varepsilon}(\delta^{x}_{\varepsilon} y ,  z))\right) \ = \ $$ 
$$= \ \mathcal{O}(\varepsilon^{2}) + d^{x}\left( \Sigma^{x}_{\varepsilon}(y,z) \,
, \, \Sigma^{x}_{\varepsilon^{2}} (y \, , \, 
\Delta^{x}_{\varepsilon}(\delta^{x}_{\varepsilon} y ,  z))\right) \quad .$$
The dilatation structure satisfies A4, therefore as $\varepsilon$ goes to $0$ we
obtain: 
$$\lim_{\varepsilon \rightarrow 0} \frac{1}{\varepsilon^{2}} \, 
Lin(x,\delta_{\varepsilon}^{x} y, \delta_{\varepsilon}^{x} z ; \varepsilon,
\varepsilon) \ = \ d^{x}\left( \Sigma^{x}(y,z) \,
, \, \Sigma^{x} (y \, , \, 
\Delta^{x}(x ,  z))\right) \ =$$ 
$$= \  d^{x}\left( \Sigma^{x}(y,z) \,
, \, \Sigma^{x}(y , z)\right) \ =  \ 0 \quad . \quad \quad \quad \square$$

The linearity of translations  
$\displaystyle \Sigma^{x}_{\varepsilon}$  is related to the linearity of the
dilatation structure, as described in the theorem below, point (a).  
As a consequence, we prove at point (b) that a linear and strong
dilatation structure comes from a conical group.

\begin{thm}
Let $(X,d,\delta)$ be a    dilatation structure. 
\begin{enumerate}
\item[(a)] If the dilatation structure is linear  then all    transformations 
$\displaystyle \Delta^{x}_{\varepsilon}(u, \cdot)$ are linear for any $u \in X$. 
\item[(b)]   If the dilatation structure is strong (satisfies A4) then it is 
linear  if and only if the dilatations  come  from the  dilatation structure of a 
 conical group, precisely for any $x \in X$ there is an open neighbourhood $D \subset X$
 of $x$ such that  $\displaystyle (\overline{D}, d^{x}, \delta)$ is the same 
 dilatation structure as the dilatation structure of the tangent space 
 of $(X,d,\delta)$ at $x$.
\end{enumerate}
\label{tdilatlin}
\end{thm}

\paragraph{Proof.}
(a) If dilatations are linear, then let $\varepsilon, \mu \in \Gamma$, $\nu(\varepsilon), 
\nu(\mu) \leq 1$, and $x, y, u, v \in X$  such that the following computations make sense. We have: 
$$\Delta^{x}_{\varepsilon} (u, \delta_{\mu}^{y} v ) = \delta_{\varepsilon^{-1}}^{\delta^{x}_{\varepsilon} u} \delta_{\varepsilon}^{x} \delta_{\mu}^{y} v \quad . $$
Let $\displaystyle A_{\varepsilon} = \delta_{\varepsilon^{-1}}^{\delta^{x}_{\varepsilon} u}$. We compute: 
$$\delta_{\mu}^{\Delta_{\varepsilon}^{x}(u,y)} \Delta_{\varepsilon}^{x} (u,v)  = \delta_{\mu}^{A_{\varepsilon} \delta_{\varepsilon}^{x} y} A_{\varepsilon} \delta_{\varepsilon}^{x} v \quad . $$
We use twice the linearity of dilatations:  
$$ \delta_{\mu}^{\Delta_{\varepsilon}^{x}(u,y)} \Delta_{\varepsilon}^{x} (u,v)  = A_{\varepsilon} 
\delta_{\mu}^{\delta_{\varepsilon}^{x} y} \delta_{\varepsilon}^{x} v = 
 \delta_{\varepsilon^{-1}}^{\delta^{x}_{\varepsilon} u} \delta_{\varepsilon}^{x} \delta_{\mu}^{y} v \quad . $$
We proved that: 
$$\Delta^{x}_{\varepsilon} (u, \delta_{\mu}^{y} v ) =  \delta_{\mu}^{\Delta_{\varepsilon}^{x}(u,y)} \Delta_{\varepsilon}^{x} (u,v)  \quad , $$
which is the conclusion of   the part (a).

(b) Suppose that the dilatation structure is strong. If dilatations are linear, then by point (a) the  transformations $\displaystyle \Delta^{x}_{\varepsilon}(u, \cdot)$$\delta$ are linear for any $u \in X$. Then, with notations made before, for $y = u$ we get 
$$\Delta^{x}_{\varepsilon} (u, \delta_{\mu}^{u} v ) =  \delta_{\mu}^{\delta_{\varepsilon}^{x} u} \Delta_{\varepsilon}^{x} (u,v)  \quad , $$
which implies 
$$\delta_{\mu}^{u} v = \Sigma^{x}_{\varepsilon} ( u, \delta^{x}_{\mu} \Delta_{\varepsilon}^{x}(u,v)) \quad . $$
We pass to the limit with $\varepsilon \rightarrow 0$ and we obtain: 
$$ \delta_{\mu}^{u} v = \Sigma^{x}( u, \delta^{x}_{\mu} \Delta^{x}(u,v)) \quad . $$
We recognize at the  right hand side the dilatations associated to the conical group 
$\displaystyle T_{x} X$. 

By proposition \ref{pexlin} the opposite implication is straightforward, 
because the dilatation structure of any conical group is linear. 
\quad $\square$

\end{document}